\documentclass[aos]{imsart}
\usepackage{booktabs}
\RequirePackage{amsthm,amsmath,amsfonts,amssymb}
\RequirePackage[authoryear]{natbib}
\RequirePackage[colorlinks,citecolor=blue,urlcolor=blue]{hyperref}
\RequirePackage{graphicx}
\usepackage{algorithmic} 
\usepackage{threeparttable}
\usepackage{bbm}
\startlocaldefs
\theoremstyle{plain}

\newtheorem{theorem}{Theorem}[section]
\newtheorem{lemma}[theorem]{Lemma}
\newtheorem{proposition}[theorem]{Proposition}

\theoremstyle{definition}

\newtheorem{remark}{Remark}

\endlocaldefs
\begin{document}

\begin{frontmatter}
\title{Learning from Uncertainty-dependent Missing Labels for Semi-supervised Classification}
\runtitle{Learning in Missing Data Analysis}
\runauthor{Y.-G. Wang et al.}

\begin{aug}
 \author[A]{\fnms{You-Gan}~\snm{Wang}
 \ead[label=e1]{you-gan.wang@uq.edu.au}$^*$\orcid{0000-0003-0901-4671}}

\author[A]{\fnms{Jinran} \snm{Wu} \ead[label=e2]{jinran.wu@uq.edu.au}\orcid{0000-0002-2388-3614}}

\author[A]{\fnms{Geoffrey J.} \snm{McLachlan} \ead[label=e3]{g.mclachlan@uq.edu.au}\orcid{0000-0002-5921-3145}}


\address[A]{School of Mathematics and Physics, The University of Queensland  \\    \printead[presep={\ }]{e1}   \printead[presep={\ }]{e2}    \printead[presep={\ }]{e3}}

\end{aug}

\begin{abstract}
Missing labels are usually regarded as a source of information loss in
classification. We study a semi-supervised setting in which the probability of
label missingness depends on the observed features through posterior
classification uncertainty. In this setting, the missingness indicator is not
only a record of an unobserved label, but also an observable signal generated by a mechanism linked to the classifier. We develop a likelihood-based information theory for such uncertainty-dependent missing labels. Under correct specification, we derive a Fisher-information
decomposition that separates a partial-labeling component from a nonnegative
mechanism-curvature term. Under joint misspecification of the label model and
the missingness mechanism, we obtain the corresponding
Godambe--Eicker--Huber--White sensitivity and sandwich-covariance partitions.
We also clarify the relevant complete-data benchmark: favorable missingness can increase information relative to ordinary fully labeled or budget-matched
non-informative labeling baselines, but cannot exceed the information in the
augmented experiment in which labels and mechanism indicators are both observed.  For plug-in classifiers, we connect the information decomposition to margin-based excess-risk bounds. In regular two-component mixture settings this yields the parametric \(n^{-1}\) excess-risk rate, with constants determined by the nuisance-adjusted information in discriminant directions. Gaussian-mixture calculations and a medical diagnosis example illustrate how
uncertainty-dependent labeling mechanisms can improve estimation and
classification under a fixed labeling budget.
\end{abstract}

\begin{keyword}[class=MSC]
\kwd[Primary ]{62F12}
\kwd{62H30}
\kwd[; secondary ]{62G20}
\kwd{68T05}
\end{keyword}

\begin{keyword}
\kwd{Semi-supervised learning}
\kwd{missing labels}
\kwd{Fisher information}
\kwd{Godambe information}
\kwd{classification}
\kwd{mixture models}
\end{keyword}

\end{frontmatter}


\maketitle
\section{Introduction}
\label{sec:intro}

Incomplete labels arise routinely in classification because expert annotation may be expensive, delayed, or feasible only for a subset of units. Missingness is therefore a familiar feature of applied data analysis, where it can reduce efficiency and complicate inference if not properly accounted for \citep[e.g.][]{Rubin1976, schafer1997analysis, LittleRubin2019}. 
Classical information arguments imply that observing only a reduced version of a fully observed experiment cannot increase information beyond the corresponding
augmented full-data experiment. This principle is fundamental, but it leaves open a more subtle possibility: when the missingness pattern is itself generated by an uncertainty-dependent mechanism, the missingness indicator may carry information about the classifier.

This possibility is central to semi-supervised classification. In many modern applications, features are inexpensive to collect whereas labels are costly, delayed, or selectively obtained. Over the past decade, semi-supervised learning (SSL) has become a central theme in statistical machine learning, driven by algorithmic ideas such as entropy regularization, consistency training, and pseudo-labeling \citep{GrandvaletBengio2004, TarvainenValpola2017, BerthelotEtAl2019, SohnEtAl2020, XieEtAl2020}, and consolidated in recent surveys \citep{yang2022survey}. In parallel, active learning emphasizes that labels are often acquired selectively, for example by querying cases for which the current classifier is most uncertain \citep{LewisGale1994, Settles2009}. In medical imaging, large annotation studies, online moderation, and decision support, triage rules and uncertainty-based heuristics may determine which units receive expensive labels and which remain unlabeled. The missingness pattern can therefore encode information about classification uncertainty.

Recent work by \citet{AhfockMcLachlan2020,ahfock2023semi} provided a key step towards an information-theoretic understanding of this phenomenon in semi-supervised mixture models. They showed that, in certain Gaussian discriminant settings, a partially labeled experiment with entropy-dependent missingness can exhibit greater observed curvature in discriminant directions than natural baseline likelihoods that do not use the missingness mechanism. This apparent paradox is resolved once the labeling mechanism is treated as part of the experiment: the missingness indicator is not merely a record of absent labels, but an observable signal generated by a mechanism tied to posterior classification uncertainty.

This paper develops a likelihood-based framework for SSL under uncertainty-dependent missingness. We study partially labeled experiments in which the probability of observing a label depends on the observed features and, through posterior class probabilities, on the label model itself. Our focus is on clarifying how such mechanisms reshape the curvature of the observed-data likelihood, when they can yield \emph{favorable missingness} relative to ordinary fully labeled or budget-matched baselines, and why this does not contradict classical information inequalities for the augmented experiment.

The perspective also connects with broader developments in machine learning. Data-collection and labeling mechanisms can create selection effects that interact with representation learning, robustness, and generalization. This view aligns with recent work linking causal reasoning and machine learning \citep{cui2020causal, scholkopf2022causality}, where selection and intervention mechanisms are treated as first-class objects rather than background nuisances. In SSL, the central challenge is therefore not merely to \emph{use} unlabeled data, but to understand when and how an uncertainty-dependent missingness process can be exploited for efficient inference and prediction, while remaining robust to model and mechanism misspecification.

The present paper extends \citet{AhfockMcLachlan2020, ahfock2023semi} in directions that are important for practice. First, we allow uncertainty-driven mechanisms, such as entropy-based masking, and joint misspecification of both the label model and the labeling mechanism, treating the likelihood as a working model. Secondly, we connect the resulting information geometry to classification performance via excess-risk bounds for plug-in classifiers fitted under uncertainty-dependent missingness.

Our main contributions are as follows. First, working directly with the observed partially labeled likelihood, we derive an information decomposition for correctly specified models that separates a partial-labeling component from a nonnegative mechanism-curvature term. Secondly, under joint misspecification, we obtain an explicit sandwich (Godambe--Eicker--Huber--White) partition for the joint estimation of label and mechanism parameters \citep{white1982maximum}, clarifying how uncertainty-dependent masking reallocates curvature in discriminant directions. Thirdly, on the prediction side, we establish margin-based excess-risk rates for plug-in classifiers, combining the asymptotic normality of the estimators with standard margin arguments \citep{Bartlett2006Tsybakov, AudibertTsybakov2007}. In regular two-component mixture settings, this yields the parametric \(n^{-1}\) excess-risk rate, with constants reflecting the contribution of the mechanism curvature. Finally, we illustrate the theory through a Gaussian mixture example and an endoscopic diagnosis study, where entropy-aware semi-supervised models achieve higher information and lower misclassification error than supervised and naive semi-supervised alternatives under comparable labeling budgets.

\section{Semi-supervised likelihood with uncertainty-dependent missingness}
\label{sec:setup}

We consider classification with \(g\) classes. Let \(Y\in\mathcal Y\subseteq\mathbb R^d\) denote the feature vector and \(Z\in\{1,\ldots,g\}\) its class label. The features are always observed, but the label may be missing. We write \(M=1\) if the label is missing and \(M=0\) if it is observed. Thus the observed data for unit \(j\) are \((Y_j,M_j,Z_j^{\rm obs})\), where \(Z_j^{\rm obs}=Z_j\) when \(M_j=0\) and $Z_j^{\rm obs}$ is unavailable otherwise.

Let \(p_\theta(Y,Z)\) denote the joint label--feature model and let
\[
  r_{\theta,\xi}(Y)=\Pr(M=1\mid Y;\theta,\xi)
\]
be the label-missingness mechanism. We allow \(r_{\theta,\xi}\) to depend on \(Y\) through posterior class probabilities
\[
  \tau_i(Y;\theta)=\Pr(Z=i\mid Y;\theta), \qquad i=1,\ldots,g,
\]
and hence on the label model itself. This is the key feature of the setting: the missingness indicator is not merely a record of absent labels, but can carry information about posterior classification uncertainty.

A useful class of mechanisms is
\[
  r_{\theta,\xi}(Y)=h\{\eta_{\theta,\xi}(Y)\},
\]
where \(h\) is a known inverse link function and
\(\eta_{\theta,\xi}(Y)\) is a predictor depending on \(\theta\), \(\xi\),
and a classification-difficulty summary of the posterior probabilities.

Let
\[
  u(Y;\theta)
\]
denote such a summary. We consider logistic-type mechanisms of the form
\[
  \eta_{\theta,\xi}(Y)
  =
  \alpha_0+\alpha_1 u(Y;\theta).
\]
Examples include the posterior Shannon entropy
\[
  e(Y;\theta)
  =
  -\sum_{i=1}^g \tau_i(Y;\theta)\log\tau_i(Y;\theta),
\]
the regularized negative log-entropy
\[
  -\log\{e(Y;\theta)+\varepsilon\},
\]
where \(\varepsilon>0\) is a small fixed constant introduced to avoid the
singularity when the posterior entropy is close to zero, and, in the binary
case, the posterior variance
\[
  \tau_1(Y;\theta)\{1-\tau_1(Y;\theta)\}.
\]
For multiclass problems, a closely related variance-type measure is the Gini
uncertainty
\[
  1-\sum_{i=1}^g \tau_i(Y;\theta)^2.
\]

The intercept \(\alpha_0\) controls the
expected missingness rate \(\gamma=E(M)\), while \(\alpha_1\) controls how
strongly missingness depends on the chosen uncertainty summary. Depending on
the definition and orientation of \(u(Y;\theta)\), the sign of \(\alpha_1\)
determines whether missingness is concentrated near high-uncertainty or
low-uncertainty cases.

The observed-data likelihood is
\begin{equation}
L(\theta,\xi)
=
\prod_{j=1}^n
\{p_\theta(Y_j,Z_j)[1-r_{\theta,\xi}(Y_j)]\}^{1-M_j}
\{p_\theta(Y_j)r_{\theta,\xi}(Y_j)\}^{M_j},
\label{eq:obs:L}
\end{equation}
where
\[
  p_\theta(Y_j)=\sum_{z=1}^g p_\theta(Y_j,z).
\]
The first factor in \eqref{eq:obs:L} is used when the label is observed; the second uses the marginal feature density when the label is missing. Estimators \((\hat\theta,\hat\xi)\) are obtained by maximizing \eqref{eq:obs:L}, or more generally by solving the corresponding estimating equations.

The mechanism \(r_{\theta,\xi}\) plays two roles. It governs the probability that a label is missing, and, when it depends on \(\theta\), it also contributes curvature to the likelihood through the observed missingness indicators \(M_j\). We call such missingness \emph{favorable} when this mechanism-induced curvature improves information in classification-relevant directions relative to a non-informative labeling scheme with the same expected labeling budget. This is a budget-matched comparison, not a claim that missingness dominates the augmented experiment in which all labels and mechanism indicators are observed.

\section{Information carried by the missingness mechanism}
\label{sec:theory}

\subsection{Ahfock–McLachlan information decomposition}

We first consider the correctly specified case. Let
\[
  S_c(Y,Z;\theta)=\nabla_\theta\log p_\theta(Y,Z)
\]
be the complete-data score, and let
\[
  S_y(Y;\theta)=E_\theta\{S_c(Y,Z;\theta)\mid Y\}
\]
be the feature-only score. Define
\[
  \mathcal I_{CC}(\theta)=E_\theta\{S_cS_c^\top\},
  \qquad
  \mathcal I_{UC}(\theta)=E_\theta\{S_yS_y^\top\},
\]
and
\[
  \mathcal I_{CC}^{({\rm miss})}(\theta)
  =
  E_\theta\!\left[
    \operatorname{Var}_\theta\{S_c(Y,Z;\theta)\mid Y\}
  \right].
\]
By the law of total covariance,
\[
  \mathcal I_{CC}(\theta)
  =
  \mathcal I_{UC}(\theta)
  +
  \mathcal I_{CC}^{({\rm miss})}(\theta).
\]

Assume that the label-missingness mechanism is MAR in the Rubin sense,
\[
  M\perp Z\mid Y,
\]
but may depend on \(\theta\) through
\[
  r_{\theta,\xi}(Y)=h\{\eta_{\theta,\xi}(Y)\}.
\]
For known \(\xi\), the observed-data Fisher information for \(\theta\) admits the decomposition
\begin{equation}
\begin{aligned}
  \mathcal I_{\rm obs}(\theta)
  &=
  \underbrace{
    \mathcal I_{CC}(\theta)
    -
    \mathcal I_{CC}^{({\rm miss},r)}(\theta,\xi)
  }_{\text{partial-labeling component}}
  +
  \underbrace{
    \mathcal I_{\rm mech}(\theta,\xi)
  }_{\text{mechanism curvature}},
\end{aligned}
\label{eq:main-decomp}
\end{equation}
where
\begin{equation}
  \mathcal I_{CC}^{({\rm miss},r)}(\theta,\xi)
  =
  E_\theta\!\left[
    r_{\theta,\xi}(Y)
    \operatorname{Var}_\theta\{S_c(Y,Z;\theta)\mid Y\}
  \right]
\label{eq:weighted-missing-info}
\end{equation}
is the missing-information loss weighted by the probability that the label is missing.

Let
\[
  \dot\eta_{\theta,\xi}(Y)=\nabla_\theta\eta_{\theta,\xi}(Y).
\]
For the logistic missingness link,
\begin{equation}
  \mathcal I_{\rm mech}(\theta,\xi)
  =
  E_\theta\!\left[
    r_{\theta,\xi}(Y)\{1-r_{\theta,\xi}(Y)\}
    \dot\eta_{\theta,\xi}(Y)\dot\eta_{\theta,\xi}(Y)^\top
  \right].
\label{eq:mech-curvature}
\end{equation}
This positive semidefinite term quantifies the curvature contributed by the observed missingness indicators.

Equivalently, the partial-labeling component in \eqref{eq:main-decomp} can be written as
\begin{equation}
\begin{aligned}
  \mathcal I_{\rm part}(\theta,\xi)
  &:=
  \mathcal I_{CC}(\theta)
  -
  \mathcal I_{CC}^{({\rm miss},r)}(\theta,\xi)  \\
  &=
  \mathcal I_{UC}(\theta)
  +
  E_\theta\!\left[
    \{1-r_{\theta,\xi}(Y)\}
    \operatorname{Var}_\theta\{S_c(Y,Z;\theta)\mid Y\}
  \right].
\end{aligned}
\label{eq:partial-equivalent}
\end{equation}
When \(r_{\theta,\xi}(Y)\equiv\gamma\), this reduces to
\[
  \mathcal I_{\rm part}(\theta,\xi)
  =
  \mathcal I_{UC}(\theta)
  +
  (1-\gamma)\mathcal I_{CC}^{({\rm miss})}(\theta).
\]

We first state and prove the information decomposition used in the main text. The key tool is Louis's observed-information identity, applied to the complete data $(Y,Z,M)$, with $(Y,M,Z^{\mathrm{obs}})$ as the observed data, where $Z^{\mathrm{obs}}=Z$ if $M=0$ and is undefined if $M=1$.

Let
\[
  S_c(Y,Z;\theta)
  =
  \nabla_\theta\log p_\theta(Y,Z),
  \qquad
  S_y(Y;\theta)
  =
  E_\theta\{S_c(Y,Z;\theta)\mid Y\}
\]
denote the complete-data and feature-only scores. Define
\[
  \mathcal I_{CC}(\theta)
  =
  E_\theta\{S_cS_c^\top\},
  \qquad
  \mathcal I_{UC}(\theta)
  =
  E_\theta\{S_yS_y^\top\},
\]
and
\[
  \mathcal I_{CC}^{(\mathrm{miss})}(\theta)
  =
  E_\theta\!\left[
    \operatorname{Var}_\theta\{S_c(Y,Z;\theta)\mid Y\}
  \right].
\]
Then, by the law of total covariance,
\[
  \mathcal I_{CC}(\theta)
  =
  \mathcal I_{UC}(\theta)
  +
  \mathcal I_{CC}^{(\mathrm{miss})}(\theta).
\]

For a feature-dependent missingness mechanism, define the weighted missing-information
loss
\[
  \mathcal I_{CC}^{(\mathrm{miss},r)}(\theta,\boldsymbol\xi)
  =
  E_\theta\!\left[
    r_{\theta,\boldsymbol\xi}(Y)
    \operatorname{Var}_\theta\{S_c(Y,Z;\theta)\mid Y\}
  \right].
\]
Also write
\[
  \dot\eta_{\theta,\boldsymbol\xi}(Y)
  =
  \nabla_\theta\eta_{\theta,\boldsymbol\xi}(Y),
  \qquad
  w_h(\eta)
  =
  \frac{\{h'(\eta)\}^2}{h(\eta)\{1-h(\eta)\}}.
\]
The corresponding mechanism-curvature term is
\[
  \mathcal I_{\mathrm{mech}}(\theta,\boldsymbol\xi)
  =
  E_\theta\!\left[
    w_h\{\eta_{\theta,\boldsymbol\xi}(Y)\}
    \dot\eta_{\theta,\boldsymbol\xi}(Y)
    \dot\eta_{\theta,\boldsymbol\xi}(Y)^\top
  \right].
\]

\begin{theorem}
\label{thm:AF}
Let \((Y,Z)\) have joint density \(p_\theta(Y,Z)\), and let
\(M\in\{0,1\}\) denote the label-missingness indicator, with \(M=1\) if the
label is missing and \(M=0\) if it is observed. Suppose that
\[
  M\perp Z\mid Y
\]
and that the missingness mechanism has the form
\[
  r_{\theta,\boldsymbol\xi}(Y)
  =
  \Pr(M=1\mid Y;\theta,\boldsymbol\xi)
  =
  h\{\eta_{\theta,\boldsymbol\xi}(Y)\},
\]
where \(h\) is a smooth inverse link function and
\(\eta_{\theta,\boldsymbol\xi}(Y)\) may depend on \(\theta\). Then, for known
\(\boldsymbol\xi\), the observed-data Fisher information for \(\theta\) satisfies
\begin{equation}
  \mathcal I_{\mathrm{obs}}(\theta)
  =
  \mathcal I_{CC}(\theta)
  -
  \mathcal I_{CC}^{(\mathrm{miss},r)}(\theta,\boldsymbol\xi)
  +
  \mathcal I_{\mathrm{mech}}(\theta,\boldsymbol\xi).
\label{eq:S1-AF-decomp}
\end{equation}
Equivalently,
\begin{equation}
  \mathcal I_{\mathrm{obs}}(\theta)
  =
  \mathcal I_{UC}(\theta)
  +
  E_\theta\!\left[
    \{1-r_{\theta,\boldsymbol\xi}(Y)\}
    \operatorname{Var}_\theta\{S_c(Y,Z;\theta)\mid Y\}
  \right]
  +
  \mathcal I_{\mathrm{mech}}(\theta,\boldsymbol\xi).
\label{eq:S1-AF-equivalent}
\end{equation}
The mechanism-curvature term is positive semidefinite. For the logistic inverse
link,
\[
  w_h\{\eta_{\theta,\boldsymbol\xi}(Y)\}
  =
  r_{\theta,\boldsymbol\xi}(Y)
  \{1-r_{\theta,\boldsymbol\xi}(Y)\}.
\]
If \(r_{\theta,\boldsymbol\xi}(Y)\equiv\gamma\), then
\[
  \mathcal I_{\mathrm{obs}}(\theta)
  =
  \mathcal I_{CC}(\theta)
  -
  \gamma\,\mathcal I_{CC}^{(\mathrm{miss})}(\theta)
  +
  \mathcal I_{\mathrm{mech}}(\theta,\boldsymbol\xi).
\]
\end{theorem}

The proof is given in Section~S1 of the Supplementary Material.  The decomposition in \eqref{eq:main-decomp} gives the Ahfock--McLachlan
information decomposition in the present entropy-driven missingness setting. 

When the mechanism is driven by posterior uncertainty, \(\mathcal I_{\rm mech}\) can be large in directions that move the decision boundary. Under a fixed labeling budget, uncertainty-dependent missingness can therefore outperform non-informative labeling by placing information where it is most useful for discrimination. In this sense, the missingness indicator acts as an observable surrogate for uncertainty about the latent label.

\subsection{Misspecification and nuisance adjustment}
\label{subsec:misspecified}

In applications, both the label model \(p_\theta(Y,Z)\) and the missingness mechanism \(r_{\theta,\xi}(Y)\) may be misspecified. We therefore treat \eqref{eq:obs:L} as a working likelihood and use the Godambe--Eicker--Huber--White sandwich framework.

We next state the  misspecification result used in the main text. Throughout this section, the observed data are
\[
        O=(Y,M,Z^{\mathrm{obs}}),
        \qquad
        Z^{\mathrm{obs}}=Z \ \text{if } M=0,
\]
and the latent full data are $(Y,Z,M)$. The observed-data score for $L_{\mathrm{obs}}(\theta,\boldsymbol\xi)$ is the conditional expectation of the latent joint score given $O$. 
Thus the Godambe partition derived below applies to the observed partially
labeled experiment. The resulting nuisance-adjusted form is the parametric analogue of the
efficient-score projection used in semiparametric missing-data theory, where
information for the target parameter is obtained after projecting away nuisance
variation in the observation mechanism
\citep{RobinsRotnitzkyZhao1994,RobinsHsiehNewey1995,Tsiatis2006}.

We now consider possible misspecification of both the label model and the
missingness mechanism. Let
\[
  O=(Y,M,Z^{\mathrm{obs}}),
  \qquad
  Z^{\mathrm{obs}}=Z \ \text{if } M=0,
\]
and let \(\vartheta=(\theta^\top,\boldsymbol\xi^\top)^\top\). Suppose the true
law \(P_0\) satisfies
\[
  p_0(y,z,m)
  =
  p_0(y,z)\,r_0(y)^m\{1-r_0(y)\}^{1-m},
  \qquad
  M\perp Z\mid Y \quad \text{under }P_0 .
\]
We fit the working joint model
\[
  p_\theta(y,z)\,
  r_{\theta,\boldsymbol\xi}(y)^m
  \{1-r_{\theta,\boldsymbol\xi}(y)\}^{1-m},
  \qquad
  r_{\theta,\boldsymbol\xi}(y)
  =
  \tilde h^{-1}\{\eta_{\theta,\boldsymbol\xi}(y)\}.
\]

Let
\[
  S_c(Y,Z;\theta)
  =
  \nabla_\theta\log p_\theta(Y,Z),
  \qquad
  S_y(Y;\theta)
  =
  E_\theta\{S_c(Y,Z;\theta)\mid Y\}
  =
  \nabla_\theta\log p_\theta(Y)
\]
be the complete-data and feature-only scores. For a general inverse link
\(\tilde h^{-1}\), define
\[
  q_{\theta,\boldsymbol\xi}(M,Y)
  =
  \frac{M-r_{\theta,\boldsymbol\xi}(Y)}
       {r_{\theta,\boldsymbol\xi}(Y)\{1-r_{\theta,\boldsymbol\xi}(Y)\}}
  \{\tilde h^{-1}\}'\{\eta_{\theta,\boldsymbol\xi}(Y)\}.
\]
Then the observed-data estimating functions are
\[
\begin{aligned}
  \psi_\theta(O;\vartheta)
  &=
  (1-M)S_c(Y,Z;\theta)
  +
  M S_y(Y;\theta)
  +
  q_{\theta,\boldsymbol\xi}(M,Y)
  \dot\eta_{\theta,\boldsymbol\xi}(Y), \\
  \psi_{\boldsymbol\xi}(O;\vartheta)
  &=
  q_{\theta,\boldsymbol\xi}(M,Y)
  \partial_{\boldsymbol\xi}\eta_{\theta,\boldsymbol\xi}(Y),
\end{aligned}
\]
where
\[
  \dot\eta_{\theta,\boldsymbol\xi}(Y)
  =
  \nabla_\theta\eta_{\theta,\boldsymbol\xi}(Y).
\]
Write
\[
  \psi(O;\vartheta)
  =
  \{\psi_\theta(O;\vartheta)^\top,
    \psi_{\boldsymbol\xi}(O;\vartheta)^\top\}^\top .
\]
Let \(\vartheta^\star=(\theta^{\star\top},\boldsymbol\xi^{\star\top})^\top\)
denote the pseudo-true value satisfying
\[
  E_0\{\psi(O;\vartheta^\star)\}=0,
\]
and define
\[
  A
  =
  -E_0\{\nabla_\vartheta\psi(O;\vartheta^\star)\},
  \qquad
  B
  =
  \operatorname{Var}_0\{\psi(O;\vartheta^\star)\}.
\]

For the \(\theta\theta\) block, set
\[
  A_{\theta\theta}^{\mathrm{pc}}
  =
  -E_0\!\left[
    \nabla_\theta
    \{(1-M)S_c(Y,Z;\theta^\star)
      +M S_y(Y;\theta^\star)\}
  \right],
\]
and
\[
  A_{\theta\theta}^{\mathrm{mech}}
  =
  -E_0\!\left[
    \nabla_\theta
    \{q_{\theta^\star,\boldsymbol\xi^\star}(M,Y)
      \dot\eta_{\theta^\star,\boldsymbol\xi^\star}(Y)\}
  \right].
\]
Also define
\[
\begin{aligned}
  B_{\theta\theta}
  =
  \operatorname{Var}_0\Big[
  &(1-M)S_c(Y,Z;\theta^\star)
    +M S_y(Y;\theta^\star)  \\
  &+
    q_{\theta^\star,\boldsymbol\xi^\star}(M,Y)
    \dot\eta_{\theta^\star,\boldsymbol\xi^\star}(Y)
  \Big].
\end{aligned}
\]

Here
\[
  \widetilde{\mathcal I}_{CC}
  =
  E_0\{S_cS_c^\top\},
  \qquad
  \widetilde{\mathcal I}_{UC}
  =
  E_0\{S_yS_y^\top\},
\]
\[
  \widetilde{\mathcal I}_{CC}^{(\mathrm{miss},r_0)}
  =
  E_0\!\left[
    r_0(Y)
    \operatorname{Var}_{\theta^\star}
    \{S_c(Y,Z;\theta^\star)\mid Y\}
  \right],
\]
and
\[
  \widetilde{\mathcal I}_{\mathrm{mech}}
  =
  E_0\!\left[
    w_{\tilde h}\{\eta_{\theta^\star,\boldsymbol\xi^\star}(Y)\}
    \dot\eta_{\theta^\star,\boldsymbol\xi^\star}(Y)
    \dot\eta_{\theta^\star,\boldsymbol\xi^\star}(Y)^\top
  \right],
\]
where
\[
  w_{\tilde h}(\eta)
  =
  \frac{\{(\tilde h^{-1})'(\eta)\}^2}
       {\tilde h^{-1}(\eta)\{1-\tilde h^{-1}(\eta)\}} .
\]

The corresponding variability block is
\[
\begin{aligned}
  B_{\theta\theta}
  =
  \operatorname{Var}_0\Big[
    &(1-M)S_c(Y,Z;\theta^\star)
    +
    M S_y(Y;\theta^\star) \\
    &+
    \{M-r_{\theta^\star,\xi^\star}(Y)\}
    \dot\eta_{\theta^\star,\xi^\star}(Y)
  \Big].
\end{aligned}
\]

Thus misspecification changes uncertainty quantification from inverse Fisher information to sandwich covariance, but it does not remove the basic partial-labeling plus mechanism-curvature structure. When \(\xi\) is unknown, comparisons for \(\theta\) should account for nuisance estimation. In the correctly specified case, this is obtained through the Schur complement
\[
  \mathcal I_{\theta\theta\cdot\xi}
  =
  \mathcal I_{\theta\theta}
  -
  \mathcal I_{\theta\xi}
  \mathcal I_{\xi\xi}^{-1}
  \mathcal I_{\xi\theta},
\]
and, under misspecification, through the corresponding \(\theta\theta\) block of \(A^{-1}BA^{-1}\). This nuisance-adjusted form is closely related to the efficient-score projection used in semiparametric missing-data theory \citep{RobinsRotnitzkyZhao1994, RobinsHsiehNewey1995, Tsiatis2006}.

\subsection{Information bounds and classification risk}
\label{subsec:upper-risk}

The preceding decomposition clarifies the role of the complete-data benchmark. If the benchmark is the ordinary fully labeled label--feature experiment based only on \((Y,Z)\), then an uncertainty-dependent missingness indicator \(M\) can contribute additional information about \(\theta\). In that sense, the observed partially labeled experiment \((Y,M,Z^{\mathrm{obs}})\) may have greater curvature in some discriminant directions than the likelihood based on \((Y,Z)\) alone. This is precisely the Ahfock--McLachlan favorable-missingness phenomenon.

There is no contradiction with classical information inequalities. If the complete experiment is augmented to include the mechanism-generated indicator \(M\), namely \((Y,Z,M)\), then the observed partially labeled data \((Y,M,Z^{\mathrm{obs}})\) are a coarsening of \((Y,Z,M)\). Hence they cannot contain more information than this augmented full experiment. Thus favorable missingness is a gain relative to ordinary fully labeled or budget-matched baselines that do not use the mechanism indicator, not a gain beyond the augmented experiment in which all labels and all mechanism indicators are observed.

\begin{theorem}
\label{thm:S-no-ex-ante}

Let \(\theta_0\) denote the true parameter of the augmented complete-data model
\(p_{\theta_0}(Y,Z,M)\), where \(M\) is the mechanism-generated missingness
indicator. Let
\[
  S_{\rm aug}(Y,Z,M;\theta_0)
  =
  \nabla_\theta \log p_\theta(Y,Z,M)\big|_{\theta=\theta_0}
\]
be the corresponding augmented complete-data score, and define
\[
  \mathcal I_{\rm aug}(\theta_0)
  =
  \operatorname{Var}\{S_{\rm aug}(Y,Z,M;\theta_0)\}.
\]
Let \(\mathcal G\) be any observed sigma-field obtained by coarsening the
augmented data, for example
\[
  \mathcal G=\sigma(Y,M,Z^{\rm obs},R),
\]
where \(R\) denotes any design randomization. Define the design information
\[
  \mathcal I_{\rm des}(\theta_0)
  =
  \operatorname{Var}\!\left[
    E\{S_{\rm aug}(Y,Z,M;\theta_0)\mid \mathcal G\}
  \right].
\]
Then
\[
  \mathcal I_{\rm des}(\theta_0)
  \preceq
  \mathcal I_{\rm aug}(\theta_0).
\]
Consequently, no regular estimator based on the coarsened observed experiment
can have asymptotic variance below the information bound for the augmented full
experiment, whenever the relevant information matrices are nonsingular. Equality
can occur only when the augmented complete-data score is
\(\mathcal G\)-measurable almost surely.
\end{theorem}

The proof is given in Section~S3 of the Supplementary Material.

\begin{lemma}[Margin exponent for regular two-component mixtures]
\label{lem:margin-two-component}
Let \(f_0\) and \(f_1\) be two regular densities on \(\mathbb R\) with equal
class priors \(1/2\), and define
\[
  \eta_0(y)
  =
  \frac{f_1(y)}{f_0(y)+f_1(y)},
  \qquad
  \Lambda(y)
  =
  \log\frac{f_1(y)}{f_0(y)}.
\]
Suppose there exists a point \(y_0\) such that \(f_1(y_0)=f_0(y_0)\), so that
\(\eta_0(y_0)=1/2\). Assume that \(\Lambda\) is differentiable in a
neighbourhood of \(y_0\), with
\[
  \Lambda'(y_0)\ne 0,
\]
and that the marginal density
\[
  f_Y(y)=\frac{1}{2}\{f_0(y)+f_1(y)\}
\]
is continuous at \(y_0\) with \(f_Y(y_0)>0\). Then there exist constants
\(C>0\) and \(t_0>0\) such that
\[
  \Pr\left\{
    \left|\eta_0(Y)-\frac12\right|\le t
  \right\}
  \le Ct,
  \qquad 0<t\le t_0.
\]
Thus the Tsybakov margin condition holds with exponent \(\alpha=1\).
\end{lemma}

The proof is given in Section~S4 of the Supplementary Material. Combining
Lemma~\ref{lem:margin-two-component} with the parametric rate
\(\hat\theta-\theta^\star=O_{\mathbb P}(n^{-1/2})\) yields
\[
  R(\hat\theta)-R(\theta^\star)=O_{\mathbb P}(n^{-1})
\]
for regular two-component mixture models with a nondegenerate decision
boundary.

\begin{proposition}[Excess-risk rate under semi-supervised designs]
\label{prop:S-excess-risk}

Let $\theta^\star$ denote the pseudo-true parameter and let $\hat\theta$ be an $M$-estimator satisfying
\[
\hat\theta-\theta^\star=O_{\mathbb P}(n^{-1/2}).
\]
Consider the plug-in classifier
\[
C_\theta(y)=\mathbbm{1}\{\mu_\theta(y)\ge 1/2\},
\]
where $\mu_\theta(y)=\tilde g^{-1}\{\eta_\theta(y)\}$ is a differentiable working regression function. Suppose that $\mu_{\theta^\star}(y)=\eta^\star(y)$ is the target regression function and that the Tsybakov margin condition holds:
\[
\Pr\!\left\{
\left|\eta^\star(Y)-\frac12\right|\le t
\right\}
\le
C_M t^\alpha
\]
for some $\alpha>0$, $C_M>0$, and all sufficiently small $t>0$. Assume also that, in a neighborhood of $\theta^\star$,
\[
\sup_y |\mu_\theta(y)-\mu_{\theta^\star}(y)|
\le
L\|\theta-\theta^\star\|
\]
for some finite constant $L$. Then
\[
R(\hat\theta)-R(\theta^\star)
=
O_{\mathbb P}\!\left(n^{-(1+\alpha)/2}\right).
\]
In particular, for the regular two-component mixture setting in
Lemma~\ref{lem:margin-two-component}, \(\alpha=1\), and hence
\[
R(\hat\theta)-R(\theta^\star)
=
O_{\mathbb P}(n^{-1}).
\]
\end{proposition}

The proof is given in Section~S5 of the Supplementary Material.
The missingness mechanism affects the constants in these rates through the
nuisance-adjusted sandwich covariance of \(\hat\theta\). Thus, when an
uncertainty-dependent mechanism increases information in discriminant
directions, it can reduce the asymptotic variance of the estimated decision
boundary and improve classification risk under the same labeling budget.

\section{Numerical illustration}
\label{sec:simulation}

This section illustrates the information decomposition in a controlled two-class Gaussian setting. The aims are to verify the mechanism-curvature term under posterior-uncertainty-dependent missingness, to identify regimes in which this term yields a net gain relative to non-informative labeling at the same budget, and to examine whether the information gain is reflected in classification performance.

Throughout, parameter estimation is carried out by maximizing the observed-data working likelihood in~\eqref{eq:obs:L}, specialized to the corresponding Gaussian mixture and a logistic missingness mechanism driven by a posterior classification-difficulty summary.

\subsection{Example: mixture of two normals}
\label{subsec:lda}

Louis~\citep{Louis1982} considered a simple but insightful example based on a
two-component Gaussian mixture. Suppose
\(Y\mid Z=i\sim \mathcal N(\mu_i,\sigma^2)\), \(i=1,2\), with mixing
proportions \(\pi_1\) and \(\pi_2=1-\pi_1\). In the information and efficiency
calculations below, the mixing proportion is treated as unknown and estimated
jointly with the component parameters.

For common variance \(\sigma^2\), define the discriminant function
\[
  d_\theta(y)
  =
  \log\frac{\pi_1}{\pi_2}
  +
  \frac{(y-\mu_2)^2-(y-\mu_1)^2}{2\sigma^2},
\]
where \(\theta=(\pi_1,\mu_1,\mu_2)\). The posterior class probability is
\[
  \tau_1(y;\theta)=h\{d_\theta(y)\},
\]
where \(h\) is the logistic function.

In the equal-prior symmetric case, \(\pi_1=\pi_2=1/2\),
\(\mu_1=\mu\), and \(\mu_2=-\mu\), this reduces to
\[
  d_\theta(y)=\frac{2\mu y}{\sigma^2}.
\]

Let \(u(y;\theta)\) denote a classification-difficulty summary computed from
the posterior probabilities. Examples include the posterior Shannon entropy
\[
  e(y;\theta)
  =
  -\sum_{i=1}^2 \tau_i(y;\theta)\log\tau_i(y;\theta),
\]
the negative log-entropy \(-\log\{e(y;\theta)+\varepsilon\}\), and the binary
posterior variance \(\tau_1(y;\theta)\{1-\tau_1(y;\theta)\}\).

We use the logistic missingness mechanism
\[
  r_\theta(y)
  =
  h\{\alpha_0+\alpha_1 u(y;\theta)\}.
\]
Here \(\alpha_0\) controls the expected missingness rate
\(\gamma=E\{r_\theta(Y)\}\), while \(\alpha_1\) controls how strongly
missingness depends on the chosen classification-difficulty summary. Depending
on the orientation of \(u(y;\theta)\), the sign of \(\alpha_1\) determines
whether missingness is concentrated near high-uncertainty or low-uncertainty
observations.

Differentiating the mechanism predictor gives
\[
  \dot\eta_\theta(y)
  =
  \alpha_1\nabla_\theta u(y;\theta).
\]
Thus the mechanism curvature along the discriminant direction \(\delta\) is
\[
  \mathcal I_{\mathrm{mech}}^{(\delta)}
  =
  \alpha_1^2
  E\!\left[
    r_\theta(Y)\{1-r_\theta(Y)\}
    \{\nabla_\delta u(Y;\theta)\}^2
  \right].
\]

For the entropy choice \(u(y;\theta)=e(y;\theta)\), the derivative is
\[
  \nabla_\delta u(y;\theta)
  =
  \tau_1(y;\theta)\{1-\tau_1(y;\theta)\}
  \log\frac{1-\tau_1(y;\theta)}{\tau_1(y;\theta)}
  \frac{y}{\sigma^2},
\]
in the equal-prior symmetric case. The corresponding expressions for
negative log-entropy or posterior variance are obtained by applying the chain
rule to the chosen \(u(y;\theta)\).

For each value of \(\alpha_1\), we tune \(\alpha_0\) so that the expected missingness rate \(\gamma\) is fixed. When \(\alpha_1>0\), increasing \(\alpha_1\) concentrates missingness near the decision boundary, where posterior uncertainty is high and the likelihood is most sensitive to the discriminant parameter.

By the weighted decomposition in Theorem~\ref{thm:AF}, the total observed information along \(\delta\) is
\[
\mathcal I_{\mathrm{obs}}^{(\delta)}
=
\mathcal I_{CC}^{(\delta)}
-
\mathcal I_{CC,\mathrm{miss},r}^{(\delta)}
+
\mathcal I_{\mathrm{mech}}^{(\delta)},
\]
where
\[
\mathcal I_{CC,\mathrm{miss},r}^{(\delta)}
=
E\!\left[
  r_\delta(Y)
  \tau_1(Y;\delta)\{1-\tau_1(Y;\delta)\}
  \frac{Y^2}{\sigma^4}
\right].
\]

The corresponding non-informative labeling baseline with the same expected missingness rate is
\[
\mathcal I_{\mathrm{MCAR}}^{(\delta)}
=
\mathcal I_{CC}^{(\delta)}
-
\gamma\mathcal I_{CC,\mathrm{miss}}^{(\delta)}.
\]

Hence the entropy-dependent mechanism improves on the same-budget non-informative baseline when
\[
\mathcal I_{\mathrm{mech}}^{(\delta)}
>
\mathcal I_{CC,\mathrm{miss},r}^{(\delta)}
-
\gamma\mathcal I_{CC,\mathrm{miss}}^{(\delta)}.
\]

Moderate values of \(\alpha_1\) yield the most efficient trade-off: small values produce little mechanism curvature, while very large values can saturate \(r_\theta(1-r_\theta)\) and reduce the net gain. This behavior is analogous to local optimal design, in which sampling effort is concentrated where the likelihood is most informative about the target direction.

\begin{remark}[Condition for favorable missingness in the Gaussian example] In the symmetric two-class Gaussian model with common variance \(\sigma^2\) and discriminant parameter \(\delta=\mu_1-\mu_2\), the complete-data score for
\(\delta\) has conditional variance
\[
  \operatorname{Var}\{S_c(Y,Z;\delta)\mid Y=y\}
  =
  \tau_1(y;\delta)\{1-\tau_1(y;\delta)\}
  \frac{y^2}{\sigma^4}.
\]
Hence the ordinary missing-information component along \(\delta\) is
\[
\mathcal I_{CC,\mathrm{miss}}^{(\delta)}
=
E\!\left[
  \tau_1(Y;\delta)\{1-\tau_1(Y;\delta)\}
  \frac{Y^2}{\sigma^4}
\right],
\]
whereas the entropy-dependent weighted loss is
\[
\mathcal I_{CC,\mathrm{miss},r}^{(\delta)}
=
E\!\left[
  r_\delta(Y)
  \tau_1(Y;\delta)\{1-\tau_1(Y;\delta)\}
  \frac{Y^2}{\sigma^4}
\right].
\]
For a logistic missingness mechanism driven by a posterior
classification-difficulty summary
\(r_\delta(y)=h\{\alpha_0+\alpha_1 e(y;\delta)\}\),
\[
\dot\eta_\delta(Y)
=
\alpha_1
\tau_1(Y;\delta)\{1-\tau_1(Y;\delta)\}
\log\!\frac{1-\tau_1(Y;\delta)}{\tau_1(Y;\delta)}
\frac{Y}{\sigma^2}.
\]
Therefore,
\[
\mathcal I_{\mathrm{mech}}^{(\delta)}(\alpha_1)
=
\frac{\alpha_1^2}{\sigma^4}
E\!\left[
  r_\delta(Y)\{1-r_\delta(Y)\}
  Y^2 A(Y;\delta)^2
\right],
\]
where
\[
A(Y;\delta)
=
\tau_1(Y;\delta)\{1-\tau_1(Y;\delta)\}
\log\!\frac{1-\tau_1(Y;\delta)}{\tau_1(Y;\delta)}.
\]

By Theorem~~\ref{thm:AF},
\[
\mathcal I_{\mathrm{obs}}^{(\delta)}(\alpha_1)
=
\mathcal I_{CC}^{(\delta)}
-
\mathcal I_{CC,\mathrm{miss},r}^{(\delta)}
+
\mathcal I_{\mathrm{mech}}^{(\delta)}(\alpha_1).
\]
Relative to a non-informative labeling design with the same missingness rate
\(\gamma\), the information advantage is

\[
\mathcal I_{\mathrm{obs}}^{(\delta)}(\alpha_1)
-
\mathcal I_{\mathrm{MCAR}}^{(\delta)}
=
\mathcal I_{\mathrm{mech}}^{(\delta)}(\alpha_1)
-
\left\{
  \mathcal I_{CC,\mathrm{miss},r}^{(\delta)}
  -
  \gamma\mathcal I_{CC,\mathrm{miss}}^{(\delta)}
\right\}.
\]
Thus the entropy-based design is favorable in the discriminant direction
relative to the same-budget non-informative baseline when
\[
\mathcal I_{\mathrm{mech}}^{(\delta)}(\alpha_1)
>
\mathcal I_{CC,\mathrm{miss},r}^{(\delta)}
-
\gamma\mathcal I_{CC,\mathrm{miss}}^{(\delta)}.
\]
The expectations are one-dimensional under the marginal mixture distribution of \(Y\) and can be evaluated numerically. This condition identifies the range of design sensitivities for which the mechanism curvature offsets the additional weighted missing-information loss relative to the same-budget non-informative baseline.
\end{remark}

For example, when \(\mu=0.6\), \(\sigma^2=1\), \(\gamma=0.5\), and \(\alpha_1=3\), a Monte Carlo calculation gives
\[
  \mathcal I_{CC}^{(\delta)}=0.25,
  \qquad
  \mathcal I_{CC}^{(\mathrm{miss}),(\delta)}=0.41,
  \qquad
  \mathcal I_{\mathrm{mech}}^{(\delta)}=0.09.
\]
The corresponding non-informative partial-labeling baseline is
\[
  \mathcal I_{\mathrm{MCAR}}^{(\delta)}
  =
  \mathcal I_{CC}^{(\delta)}
  -
  \gamma\mathcal I_{CC}^{(\mathrm{miss}),(\delta)}
  =
  0.04.
\]
Adding the mechanism-curvature term gives
\[
  \mathcal I_{\mathrm{MCAR}}^{(\delta)}
  +
  \mathcal I_{\mathrm{mech}}^{(\delta)}
  =
  0.13.
\]
Thus, at the same missingness rate, the uncertainty-dependent mechanism gives about a threefold increase in discriminant-direction information relative to the non-informative baseline. This gain should be interpreted in the budget-matched sense: it arises because the observed missingness indicators carry information about posterior uncertainty near the decision boundary.

We next examine this behavior over a grid of separations, missingness rates, and design sensitivities. We take \(Y\in\mathbb R\) and \(Z\in\{1,2\}\), with equal priors \(\pi_1=\pi_2=1/2\). Conditional on \(Z=i\),
\[
  Y\mid Z=i\sim\mathcal N(\mu_i,\sigma^2),
  \qquad i=1,2.
\]
The Mahalanobis separation $ \delta=\frac{\mu_1-\mu_2}{\sigma}$ controls class overlap. We use the symmetric setting \(\mu_1=+\mu\), \(\mu_2=-\mu\), and \(\sigma^2=1\), so that \(\delta=2\mu\). The posterior probability is
\[
  \tau_1(y;\mu)=\eta(2\mu y),
\]
where \(\eta(\cdot)\) is the logistic function, and entropy is
\[
  e(y;\mu)
  =
  -\sum_{i=1}^2 \tau_i(y;\mu)\log\tau_i(y;\mu).
\]
Labels are made missing according to
\[
  r_\mu(y\mid\alpha_0,\alpha_1)
  =
  \Pr(M=1\mid Y=y)
  =
  \eta\{\alpha_0+\alpha_1 e(y;\mu)\}.
\]
We explore
\[
  \delta\in\{1,2,3\},
  \qquad
  \gamma\in\{0.1,0.3,0.5\},
  \qquad
  \alpha_1\in\{2,4,8\},
\]
tuning \(\alpha_0\) to achieve the desired missingness rate.

\begin{figure}[!ht]
\centering
\includegraphics[width=\linewidth]{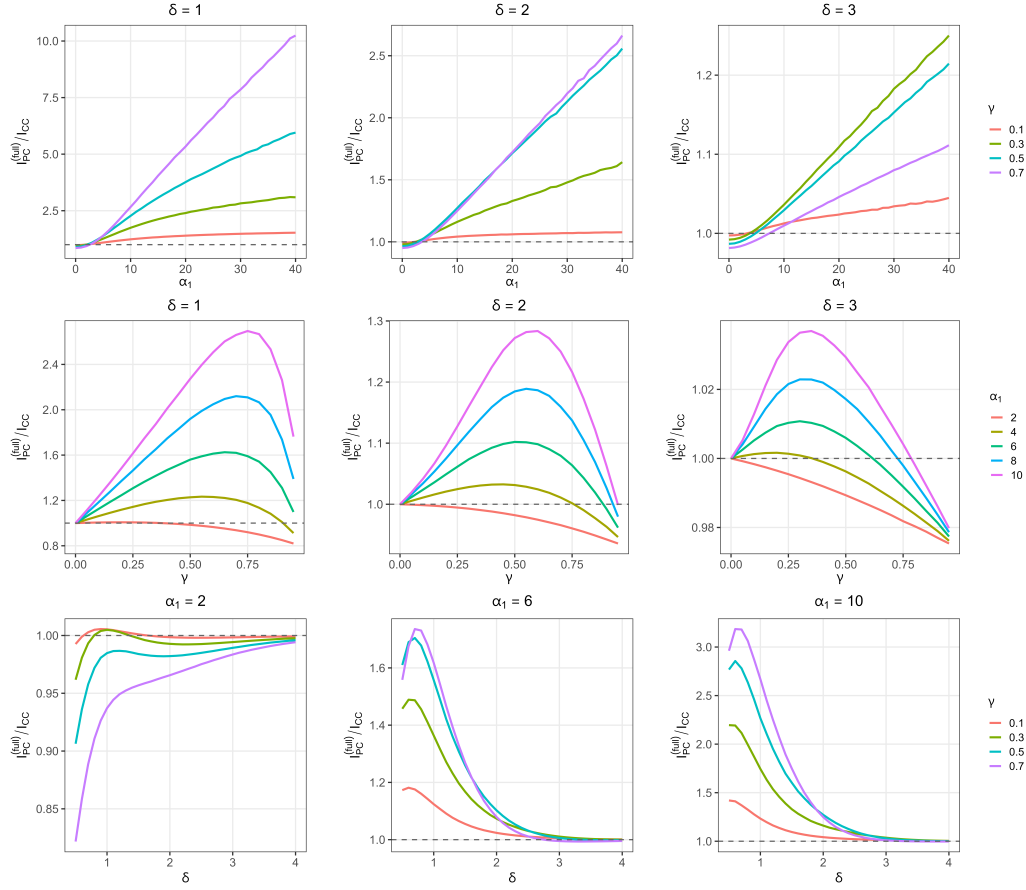}
\caption{\textbf{Information gain ratio.} Panels~A--C illustrate the behavior of the information-gain ratio \(R(\delta,\gamma,\alpha_1)\) under varying design parameters. \textbf{(A)}~\(R\) versus the design sensitivity \(\alpha_1\) for \(\delta\in\{1,2,3\}\) and \(\gamma\in\{0.1,0.3,0.5,0.7\}\). \textbf{(B)}~\(R\) versus the missingness proportion \(\gamma\) for \(\alpha_1\in\{2,4,6,8,10\}\). \textbf{(C)}~\(R\) versus the class-separation parameter \(\delta\in[0.5,4.0]\) for \(\alpha_1\in\{2,6,10\}\) and \(\gamma\in\{0.1,0.3,0.5,0.7\}\). Each curve represents an average over \(B=100\) Monte Carlo replicates with \(n=2000\) observations. The dashed horizontal line marks \(R=1\); values above this threshold indicate that the observed semi-supervised likelihood, including the missingness indicators, contains more information than the corresponding label--feature benchmark in the plotted direction.}
\label{fig:R-vs-alpha1}
\end{figure}

Figure~\ref{fig:R-vs-alpha1} reports Monte Carlo estimates of \(R(\delta,\gamma,\alpha_1)\). For moderate class separation and moderate missingness, the ratio exceeds one for sufficiently strong but non-saturated design sensitivity, indicating that the entropy-dependent mechanism contributes substantial curvature. Panel~A shows that increasing \(\alpha_1\) can improve information when the classes are not already well separated. Panel~B shows a non-monotone dependence on the missingness rate: some missingness can be informative when it is uncertainty-dependent, but excessive missingness removes too many labels. Panel~C shows that the gain is largest at intermediate class separation, where the decision boundary is neither trivial nor unidentifiable.

Overall, the simulation confirms that favorable missingness is a local and regime-dependent phenomenon. It is strongest when three conditions are balanced: moderate class overlap, a non-excessive missingness rate, and a mechanism sensitive enough to concentrate missingness near high-uncertainty observations without saturating.

\subsection{Numerical simulations}

We next adopt the two-component mixture setting of Example~2 in~\cite{Louis1982}. Data are generated from a symmetric Gaussian mixture with $\pi_1=\pi_2=0.5$, $\mu_1=2$, $\mu_2=0$, and $\sigma_1^2=\sigma_2^2=1$, with total sample size $n=500$. Entropy is computed under the true parameter values, and the missing-label indicator is generated from the logistic mechanism with $\alpha_0=-3$ and $\alpha_1=5$. This produces a partially labeled dataset consisting of feature vectors, observed labels for non-missing cases, and missingness indicators.

To compare directly with the Fisher-information decomposition, we estimate the mixing proportion $\pi_1$ and component means $\mu_1$ and $\mu_2$, holding the variances and missingness parameters fixed at their true values. The estimates are
\[
\hat{\pi}_1=0.51,\qquad
\hat{\mu}_1=1.98,\qquad
\hat{\mu}_2=0.01,
\]
with $\hat\gamma=0.286$. Table~\ref{tab:fisher-decomp-mc} reports the corresponding information components. The mechanism-curvature matrix is large enough to offset much of the missing-information loss, so that the observed likelihood has substantially higher curvature than the non-informative partial-labeling baseline.

\begin{table}[!ht]
\centering
\caption{Fisher information decomposition under entropy-based MAR missingness}
\label{tab:fisher-decomp-mc}
\renewcommand{\arraystretch}{1.25}
\begin{tabular}{c c}
\hline
Information component & Fisher information matrix \\
\hline
$\mathcal I_{CC}(\hat{\theta})$
&
$\begin{pmatrix}
2001 & 0 & 0 \\
0 & 254 & 0 \\
0 & 0 & 246
\end{pmatrix}$
\\[6pt]
$\hat{\gamma} \mathcal I^{(\mathrm{miss})}_{CC}(\hat{\theta})$
&
$\begin{pmatrix}
390 & -97 & -95 \\
-97 & 32 & 16 \\
-95 & 16 & 31
\end{pmatrix}$
\\[6pt]
$\mathcal I_{\mathrm{mech}}(\hat{\theta},\xi)$
&
$\begin{pmatrix}
1017 & -238 & -262 \\
-238 & 103 & 13 \\
-262 & 13 & 261
\end{pmatrix}$
\\[6pt]
$\mathcal I_{\mathrm{obs}}(\hat{\theta})$
&
$\begin{pmatrix}
2628 & -141 & -167 \\
-141 & 325 & -2 \\
-167 & -2 & 476
\end{pmatrix}$
\\
\hline
\end{tabular}
\end{table}

\subsection{A cost--benefit example}

We finally examine how uncertainty-dependent missingness interacts with a fixed labeling budget. Let the cost of collecting features for one unit be 1, and let the additional cost of obtaining its label be $c$. If $n$ feature vectors are collected and $n_\ell$ of them are labeled, the total cost is
\[
C=n+c n_\ell.
\]
A fully supervised design has $n=n_\ell$ and therefore sample size
\[
n_0=\frac{C}{1+c}.
\]
In contrast, a semi-supervised design can collect more feature vectors while labeling only a fraction of them. Let $\gamma$ denote the missingness rate, so that
\[
n_\ell=n(1-\gamma).
\]
Substituting this into the cost constraint gives
\[
n=\frac{C}{1+c(1-\gamma)},\qquad
n_\ell=\frac{(1-\gamma)C}{1+c(1-\gamma)}.
\]
Thus, increasing $\gamma$ reduces the number of labeled observations but allows a larger feature sample under the same total budget. The optimal value of $\gamma$ balances the loss of labels against the information carried by the uncertainty-dependent missingness mechanism.

\begin{figure}[!ht]
    \centering
    \includegraphics[width=0.75\linewidth]{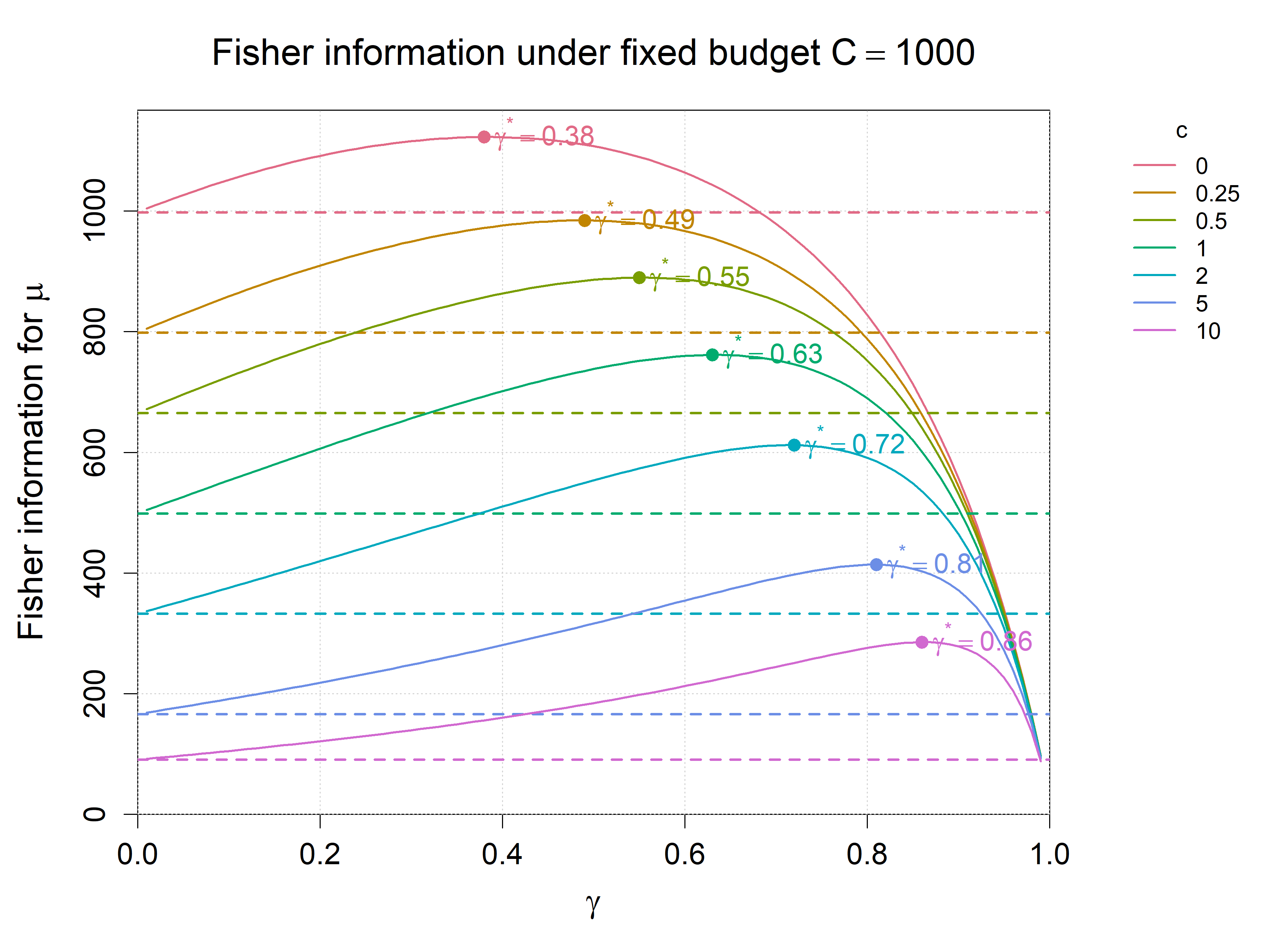}
  \caption{Nuisance-adjusted Fisher information for the mean parameter \(\mu\)
as a function of the missingness rate \(\gamma\) under a fixed total budget
\(C=1000\), for different labeling cost ratios
\(c\in\{0,0.25,0.5,1,2,5,10\}\). The mixing proportion, component means, and
variance parameters are estimated jointly. The dashed horizontal line indicates
the complete-case baseline. The optimal missing proportion under each labeling
cost ratio is denoted by \(\gamma^\ast\). The results illustrate how different
design strategies trade off information efficiency under the same cost
constraint.}
    \label{fig:budget-fi}
\end{figure}

For illustration, we set $C=1000$ and consider labeling cost ratios $c\in\{0,0.25,0.5,1,2,5,10\}$. Data are generated from a symmetric two-component Gaussian mixture with $\pi_1=\pi_2=0.5$, means $\mu_1=-\mu_2=0.5$, and equal variances $\sigma_1^2=\sigma_2^2=1$. The mechanism sensitivity is fixed at $\alpha_1=3$. 
The mixing proportion, component means, and variance parameters are estimated
jointly. Thus the Fisher information underlying Figure~\ref{fig:budget-fi} is
computed for the full Gaussian-mixture parameter vector relevant to
classification, rather than for the mean parameter alone.

Figure~\ref{fig:budget-fi} shows the resulting information, or relative
efficiency, as a function of the missingness rate \(\gamma\). The optimal
missingness rate increases with the relative cost of labeling, but its value
also reflects the additional uncertainty from estimating the class prior
probability and variance parameters. The figure therefore illustrates how
uncertainty-dependent semi-supervised designs trade off feature sample size,
labeling cost, mechanism curvature, and full mixture-parameter estimation under
a common budget constraint.

\section{A case study}
\label{sec:case-study}

We illustrate the proposed framework using the gastrointestinal dataset of~\citet{mesejo2016computer}. The dataset contains 76 colonoscopy videos, each with a histology-based ground-truth label and diagnostic assessments from seven endoscopists, including four experts and three beginners. Both white-light and narrow-band imaging modalities were used to classify each lesion as benign or malignant. For each case, four representative image-derived features, numbered 294, 441, 472, and 486, were extracted from the videos.

Each endoscopist independently judged whether the patient required resection (malignant) or no resection (benign). To construct a partially labeled sample, we treated observations with unanimous agreement among all seven endoscopists as labeled, assigning their labels according to the histology ground truth: 1 for ``resection'' and 2 for ``no resection''. Observations without unanimous agreement were treated as having missing labels. This yielded 35 labeled observations, consisting of 31 ``resection'' and 4 ``no resection'' cases. In the full histology-labeled dataset, there are 55 ``resection'' and 21 ``no resection'' cases.

To visualize the missingness pattern, we first fitted an equal-covariance Gaussian classification model to the labeled observations by minimizing the complete-data negative log-likelihood. Entropy values were then computed from the fitted posterior class probabilities, and the missingness indicator was modeled as a function of entropy using logistic, cauchit, log-log, and Gaussian-kernel specifications.

\begin{figure}[!ht]
\centering
\includegraphics[width=0.95\linewidth]{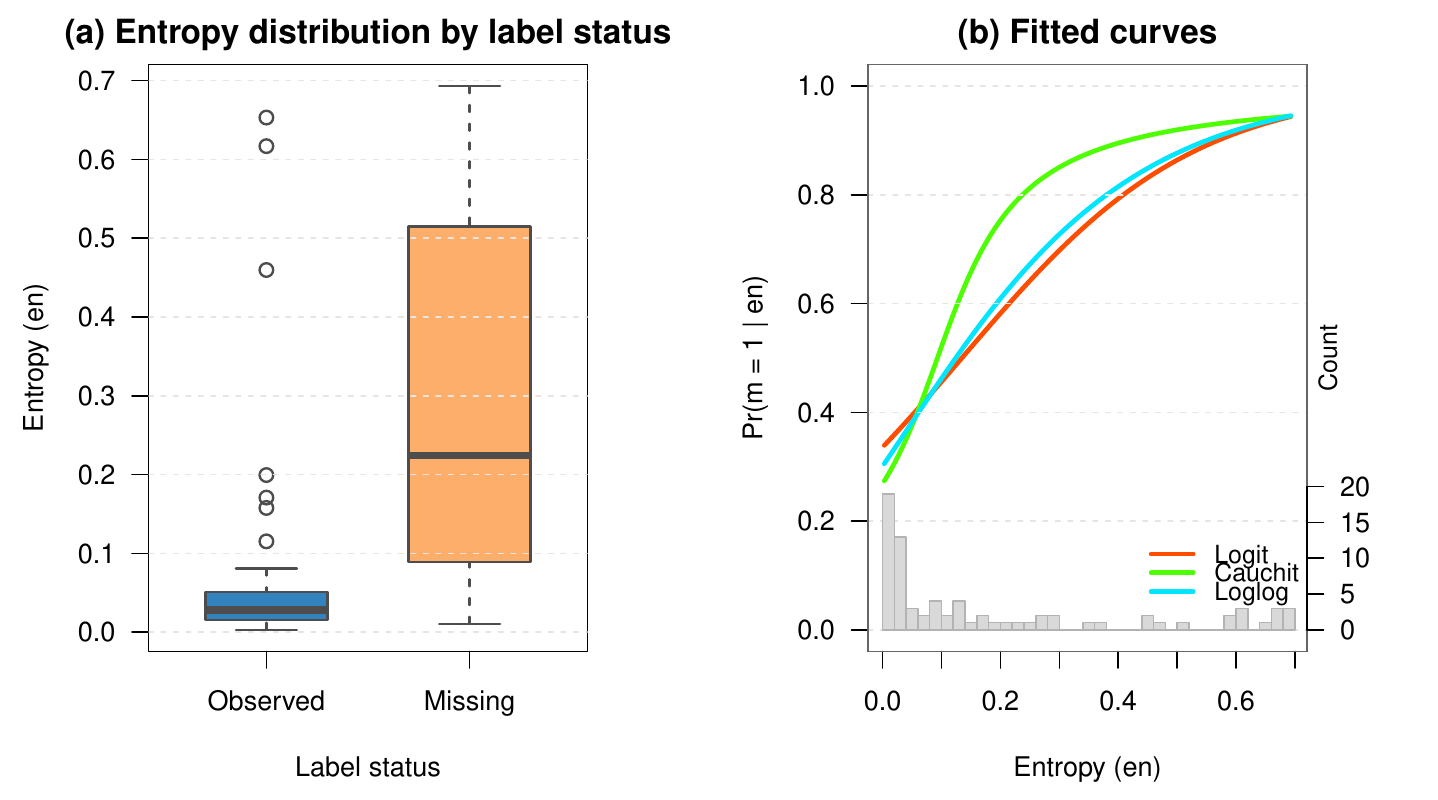}
\caption{Relationship between sample entropy and label missingness in the gastrointestinal dataset. The left panel compares the entropy distributions for observations with and without labels, while the right panel depicts the fitted curves linking the probability of missingness to entropy under different models.}
\label{missingness_pattern}
\end{figure}

Figure~\ref{missingness_pattern} shows the resulting relationship between entropy and label missingness. The left panel indicates that observations with missing labels tend to have higher entropy than observations with observed labels. The right panel shows that the estimated probability of missingness generally increases with entropy across the fitted parametric mechanisms. This supports the working assumption that label missingness is related to posterior classification uncertainty.

\begin{table*}[!ht]
\centering
\caption{Summary of model comparison, Fisher information, and predictive error rates for the gastrointestinal dataset. Panel~(a) reports model-selection indices \((-2\mathrm{LL})\), AIC, and BIC; Panel~(b) presents parameter-wise Fisher information; Panel~(c) shows average classification error rates over 200 replications. Higher log-likelihood and lower AIC/BIC indicate better model fit, larger Fisher information indicates higher estimation efficiency, and lower error rates indicate better predictive performance.}
\label{tab:case_summary}
\begin{threeparttable}
\renewcommand{\arraystretch}{1.15}

(a) Model selection criteria \\[0.3em]
\makebox[\linewidth][c]{
\begin{tabular}{lccc}
\toprule
Model & \(-2\)LL & AIC & BIC \\
\midrule
\(\mathrm{SSL}^{\text{logit}}\)   & -461.6081 & \(-419.6081\) & \(-370.6627\) \\
\(\mathrm{SSL}^{\text{cauchit}}\) & -460.1512 & \(-418.1512\) & \(-369.2058\) \\
\(\mathrm{SSL}^{\text{loglog}}\)  & -461.3591 & \(-419.3591\) & \(-370.4137\) \\
\bottomrule
\end{tabular}%
}

\vspace{1em}

(b) Fisher information components \\[0.3em]
\makebox[\linewidth][c]{%
\resizebox{0.95\linewidth}{!}{
\begin{tabular}{lccccccccc}
\toprule
Model & \(\pi_1\) & \(\mu_{1,1}\) & \(\mu_{1,2}\) & \(\mu_{1,3}\) & \(\mu_{1,4}\) &
       \(\mu_{2,1}\) & \(\mu_{2,2}\) & \(\mu_{2,3}\) & \(\mu_{2,4}\) \\
\midrule
SL (complete data)
& 380.56 & 5197.26 & 5253.20 & 10761.87 & 8221.89 & 1984.41 & 2005.77 & 4109.08 & 3139.27 \\
\(\mathrm{SSL}^{\text{logit}}\)
& \textbf{536.91} & \textbf{6762.60} & \textbf{5636.11} & \textbf{19073.89} &
  \textbf{9563.01} & \textbf{3217.31} & \textbf{3170.42} & \textbf{7639.12} & \textbf{4126.68} \\
\(\mathrm{SSL}^{\text{cauchit}}\)
& 451.47 & 5876.33 & 5458.69 & 16634.87 & 9402.28 & 2349.62 & 2106.07 & 6661.55 & 3951.19 \\
\(\mathrm{SSL}^{\text{loglog}}\)
& 505.52 & 6365.66 & 5530.43 & 17759.23 & 9227.48 & 2842.39 & 2747.94 & 7143.67 & 3885.40 \\
\bottomrule
\end{tabular}%
}}

\vspace{1em}

(c) Prediction error rates \\[0.3em]
\makebox[\linewidth][c]{
\begin{tabular}{lcccc}
\toprule
Model & \(\mathrm{SSL}^{\text{cauchit}}\) & \(\mathrm{SSL}^{\text{logit}}\) &
\(\mathrm{SSL}^{\text{loglog}}\) & SL (complete data) \\
\midrule
Error rate & 0.1338 & \textbf{0.1325} & 0.1413 & 0.1775 \\
\bottomrule
\end{tabular}
}

\end{threeparttable}
\end{table*}

We then fitted the partially labeled dataset using the observed-data likelihood under an entropy-based MAR mechanism with equal covariance matrices. The missingness model was specified using three binomial links: logit, cauchit, and log-log. The corresponding models are denoted by \(\mathrm{SSL}^{\text{logit}}\), \(\mathrm{SSL}^{\text{cauchit}}\), and \(\mathrm{SSL}^{\text{loglog}}\), respectively.

Table~\ref{tab:case_summary} summarizes the model-selection criteria, parameter-wise Fisher information, and prediction error rates. Among the three entropy-based semi-supervised models, \(\mathrm{SSL}^{\text{logit}}\) gives the best overall fit, with the lowest \(-2\mathrm{LL}\), AIC, and BIC.

We also compared the model-based Fisher information for the mixing proportion \(\pi_1\) and the component means \(\boldsymbol\mu_1=(\mu_{1,1},\mu_{1,2},\mu_{1,3},\mu_{1,4})\) and \(\boldsymbol\mu_2=(\mu_{2,1},\mu_{2,2},\mu_{2,3},\mu_{2,4})\). Panel~(b) shows that \(\mathrm{SSL}^{\text{logit}}\) has the largest observed information for most parameters, followed by \(\mathrm{SSL}^{\text{loglog}}\) and \(\mathrm{SSL}^{\text{cauchit}}\).

The semi-supervised models also show larger fitted information than the supervised benchmark using ground-truth labels. This comparison should be interpreted as a model-based curvature comparison under the fitted working likelihood: the gain arises from the additional curvature contributed by the entropy-dependent missingness indicators, not from a violation of the augmented complete-data information bound.

To assess predictive performance, we randomly selected 90\% of the observations for training and used the remaining 10\% for testing. Classification performance was evaluated against the histology ground-truth labels, and the split was repeated 200 times. The average error rates are reported in Panel~(c) of Table~\ref{tab:case_summary}.

The supervised benchmark has the largest error rate, 0.1775, whereas all three entropy-based semi-supervised models perform better. Among them, \(\mathrm{SSL}^{\text{logit}}\) achieves the lowest average error rate, 0.1325, followed closely by \(\mathrm{SSL}^{\text{cauchit}}\) at 0.1338 and \(\mathrm{SSL}^{\text{loglog}}\) at 0.1413.

Overall, the three panels give a coherent picture. The missingness pattern is strongly associated with posterior uncertainty, the logit mechanism gives the best-fitted likelihood among the entropy-based models, and the same model also yields the largest model-based information and the lowest prediction error. The case study therefore illustrates how modeling uncertainty-dependent missingness can improve both estimation and classification under a fixed partially labeled dataset.

\section{Discussion}
\label{sec:discussion}

Semi-supervised learning (SSL) is motivated by a common asymmetry in scientific data: features are often easy to collect, whereas labels may be costly, delayed, selectively observed, or available only after expert review. Much of the SSL literature therefore develops algorithmic strategies---including entropy regularization, consistency training, and pseudo-labeling---to exploit unlabeled data for prediction \citep{GrandvaletBengio2004, TarvainenValpola2017, BerthelotEtAl2019, SohnEtAl2020, XieEtAl2020, yang2022survey}. Active learning provides a complementary perspective, emphasizing that labels are often acquired selectively, for example by querying cases for which the current classifier is most uncertain \citep{LewisGale1994, Settles2009}. This paper takes a likelihood-based view of the same phenomenon: we ask how the \emph{pattern} of missing labels, and its dependence on features or posterior classification uncertainty, changes the information available for estimation and classification.

The central point is that uncertainty-dependent missingness is not merely ordinary MAR in a passive sense. Although the mechanism is MAR in the Rubin terminology when \(M\perp Z\mid Y\), it is a structured MAR mechanism tied to the classifier itself. When the probability of observing a label depends on posterior uncertainty, the missingness indicator becomes an observable signal about where the latent class label is difficult to infer. The resulting observed-data information separates into a partial-labeling component and a nonnegative mechanism-curvature term. This term quantifies the information carried by the missingness indicators and explains the Ahfock--McLachlan favorable-missingness phenomenon: under a fixed labeling budget, uncertainty-dependent missingness can increase information in discriminant directions relative to non-informative labeling at the same budget.

Under misspecification of the label model or the missingness mechanism, uncertainty quantification is governed by the sandwich covariance in the Godambe--Eicker--Huber--White framework \citep{Godambe1960,Eicker1963, Huber1964, white1982maximum}. The corresponding sensitivity partition retains the same partial-labeling plus mechanism-curvature structure, while the sandwich covariance provides the appropriate basis for standard errors and confidence regions. This is important in applications, where entropy-based or uncertainty-based mechanisms are usually working models rather than exact descriptions of the data-generating process.

The Gaussian mixture experiments and the endoscopic diagnosis study illustrate the main message. Explicitly modeling uncertainty-dependent missingness can improve both fitted information and classification performance relative to supervised or naive semi-supervised baselines with the same labeling budget. These gains do not contradict classical information inequalities. As clarified above, favorable missingness is a comparison with ordinary fully labeled or budget-matched baselines that do not use the mechanism indicator. When the experiment is augmented to include \(M\), the observed partially labeled data are still a coarsening of \((Y,Z,M)\). Thus the mechanism does not create information \emph{ex nihilo}; rather, a genuine uncertainty-dependent mechanism makes \(M\) an additional observable signal whose dependence on \(Y\) can be exploited.

The framework also connects to the selective-labels problem, where outcomes are observed only for a non-representative subset of units \citep{LakkarajuEtAl2017}. In such settings, the observation process is not merely a nuisance but part of the statistical experiment. Explicitly modeling it can improve efficiency and interpretability, while also clarifying when an apparent information gain is due to a meaningful uncertainty-dependent mechanism rather than to model artifacts. A practical extension is to use bootstrap diagnostics to assess the stability of the curvature decomposition. Parametric or semi-parametric bootstrap procedures under the fitted joint model could provide standard errors for information components, assess confidence regions for discriminant contrasts, and compare uncertainty-dependent masking with simpler MCAR baselines at the same missingness rate \citep{Efron1979, EfronTibshirani1993, DavisonHinkley1997}.

Several limitations point to future work. The present analysis is local and information-based, with main developments for parametric mixtures and GLM-type working models. Extending the theory to modern deep SSL, where curvature must be approximated by automatic differentiation and stochastic optimization, is a natural next step. It would also be useful to embed the decomposition in a decision-analytic framework that explicitly trades off predictive risk, labeling cost, and mechanism sensitivity, thereby connecting the present likelihood-based view more closely to active learning \citep{LewisGale1994, Settles2009}. More broadly, recent work at the interface of causal inference and machine learning emphasizes that selection and intervention mechanisms can be central to robust generalization \citep{cui2020causal, scholkopf2022causality}; uncertainty-dependent labeling
provides a concrete setting in which these issues arise.

Finally, uncertainty-driven labeling raises representativeness and fairness questions. If high-uncertainty regions coincide with under-represented groups, rare subtypes, or clinically important boundary cases, selective labeling may change the effective training distribution and affect downstream decisions. Future work should therefore combine information-based diagnostics with group-wise performance checks and design constraints, so that efficiency gains from favorable missingness do not come at the expense of equity, coverage, or scientific validity.


\bibliographystyle{plainnat}
\bibliography{refs}

\end{document}